\documentclass[a4paper,11pt]{amsart}

\usepackage{amsmath,amsfonts,amssymb,amsthm,latexsym,amscd}
\usepackage[dvips]{graphicx}
\usepackage[T1]{fontenc}

\usepackage{parskip}

\newtheorem{thm}{Theorem}[section]
\newtheorem{thm*}{Theorem}
\newtheorem{lem}[thm]{Lemma}
\newtheorem{prop}[thm]{Proposition}
\newtheorem{cor}[thm]{Corollary}
\newtheorem{defn}[thm]{Definition}
\newtheorem{rem}[thm]{Remark}
\newtheorem{defn*}[thm*]{Definition}
\newtheorem{rem*}[thm*]{Remark}
\newtheorem{pr*}[thm*]{Problem}
\newtheorem{note*}[thm*]{Remark}

\def\C{\mathcal{C}}
\def\cD{\mathcal{D}}
\def\D{\mathbb{D}}
\def\cH{\mathcal{H}}
\def\Hom{\operatorname{Hom}}
\def\N{\mathbb{N}}
\def\R{{\mathbb{R}}}
\def\cS{\mathcal{S}}
\def\Z{\mathbb{Z}}
\def\a{\alpha}
\def\area{\operatorname{area}}
\def\b{\beta}
\def\g{\gamma}

\def\lk{\operatorname{lk}}
\def\o{\omega}
\def\s{\sigma}
\def\vol{\operatorname{vol}}

\begin{document}

\title{On quasi-morphisms from knot and braid invariants}
%Short title: On quasi-morphisms from knot invariants.

\author{Michael Brandenbursky}

%\address{\it{Department of Mathematics, Vanderbilt University, Nashville, TN 37240, USA}\\
%\it{michael.brandenbursky@vanderbilt.edu}}

\maketitle

\begin{abstract}
We study quasi-morphisms on the groups $P_n$ of pure braids on $n$ strings and on the group $\cD$ of compactly supported
area-preserving diffeomorphisms of an open two-dimensional disc. We show that it is possible to build quasi-morphisms on $P_n$
by using knot invariants which satisfy some special properties. In particular, we study quasi-morphisms which come from knot
Floer homology and Khovanov-type homology. We then discuss possible variations of the Gambaudo-Ghys construction, using the above quasi-morphisms on $P_n$ to build quasi-morphisms on the group $\cD$ of diffeomorphisms of a 2-disc.
\end{abstract}

%\keywords{quasi-morphisms, braid groups, knot concordance invariants, area-preserving diffeomorphisms.}

%\ccode{Mathematics Subject Classification 2000: 20F36, 37E30, 57M07, 57M25, 57M27}

\section{Introduction}

Real-valued quasi-morphisms are known to be a helpful tool in the
study of algebraic structure of non-Abelian groups, especially the
ones that admit a few or no (linearly independent) real-valued homomorphisms.
Recall that a {\it quasi-morphism} on a group $G$ is a function
$\varphi: G \to \R$ which satisfies the homomorphism equation up to
a bounded error: there exists $K_\varphi > 0$ such that
$$|\varphi (ab) -\varphi(a) -\varphi (b)| \leq K_\varphi$$ for all $a,b
\in G$. A quasi-morphism $\varphi$ is called {\it homogeneous} if
$\varphi (a^m) = m \varphi (a)$ for all $a \in G$ and $m \in \Z$.
Any quasi-morphism $\varphi$ can be {\it homogenized}: setting
\begin{equation}\label{eq:definition-homogenization}
\widetilde{\varphi} (a) := \lim_{k\to +\infty} \varphi (a^k)/k
\end{equation}
we get a homogeneous (possibly trivial) quasi-morphism
$\widetilde{\varphi}$.

In \cite{surfaces} Gambaudo and Ghys showed that a signature invariant of links in $S^3$ defines quasi-morphisms on the full braid groups $B_n$ on $n$ strings. These quasi-morphisms on $B_n$ are
constructed in the following way: close up a braid to a link in
the standard way and take the value of the signature invariant on that link.
In this paper we are motivated by the following question:"What happens if one plugs in a general knot/link invariant
in this construction?" In this paper we give a sufficient condition which is satisfied by some already known
knot invariants. More specifically we show that any homomorphism from the concordance group of knots in $S^3$ to the reals, which is bounded in some canonical norm on this group, defines a quasi-morphism on $B_n$. Known knot invariants which satisfy the above condition \emph{do not produce new quasi-morphisms} on $B_n$. However, this condition \emph{may possibly} lead to new examples of quasi-morphisms.

We consider three specific remarkable knot/link invariants: the  Rasmussen link invariant $s$ \cite{Beliakova,Rasmussen},
which comes from a Khovanov-type theory, the Ozsvath-Szabo knot invariant $\tau$ \cite{Ozsvath-tau}, which
comes from the knot Floer homology, and the classical signature link invariant $sign$ \cite{Likorish,Rolfsen}. In \cite{Baader}
Baader has shown that the Rasmussen link invariant $s$ defines a quasi-morphism on $B_n$. We show that the homogenization
of this quasi-morphism is equal to the classical linking number homomorphism $lk$ on $B_n$. We also show that the situation with the
Ozsvath-Szabo knot invariant $\tau$ is similar: it defines a quasi-morphism on $B_n$ and its homogenization is again the linking number
homomorphism divided by $2$. In addition we show that the homogenization of an induced signature quasi-morphism on $B_n$ and $lk$ coincide on alternating braids. We also present an inequality which connects $s$, $\tau$ and the braid index of a knot.

Further we discuss the group $\cD$ of compactly supported
area-preserving diffeomorphisms of the open unit disc in the Euclidean plane. The group $\cD$ admits a unique
(continuous, in the proper sense) homomorphism to the reals -- the famous Calabi homomorphism (see e.g. \cite{Banyaga,Calabi,enlacement}). At the same time $\cD$ is known to admit many (linearly independent) homogeneous quasi-morphisms
(see e.g. \cite{Barge-Ghys,Biran-Entov-Polterovich,surfaces}). In this work we consider a
particular geometric construction of such quasi-morphisms, essentially contained in \cite{surfaces}, which produces quasi-morphisms on $\cD$ from
quasi-morphisms on the pure braid groups $P_n$. We discuss the computation of the quasi-morphisms on $\cD$, obtained by this construction, on diffeomorphisms generated by time-in\-de\-pen\-dent (compactly supported) Hamiltonians. For a generic Hamiltonian $H$ of this sort we
present the result of the computation in terms of the Reeb graph of $H$ and the integral of the push-forward of $H$ to the graph
against a certain signed measure on the graph. This result enables us to show that the Calabi homomorphism and the quasi-morphism on
$\cD$ induced by the signature invariant of $n$-component links are asymptotically equivalent, as $n\to \infty$, on the flows
generated by time-independent (compactly supported) Hamiltonians.
\vspace{0.7mm}

\textbf{Plan of the paper.}
In Section 2 we  formulate sufficient conditions for a knot invariant to yield a quasi-morphism on $B_n$.
In Section 3 we give examples of such knot invariants and provide properties of homogeneous quasi-morphisms, defined by them, on $B_n$.
In Section 4 we discuss the Gambaudo-Ghys construction, which produces homogeneous quasi-morphisms on $\cD$ from homogeneous quasi-morphisms on $P_n$.
We define a set of generic autonomous Hamiltonians, and discuss the computation of the induced quasi-morphisms on $\cD$ on the elements generated by these Hamiltonians. At the end we discuss the asymptotic behavior of the induced signature quasi-morphism.
\vspace{0.7mm}

\section{Quasi-morphisms on braid groups defined by knot invariants}

\begin{figure}[htb]
\centerline{\includegraphics[height=1.7in]{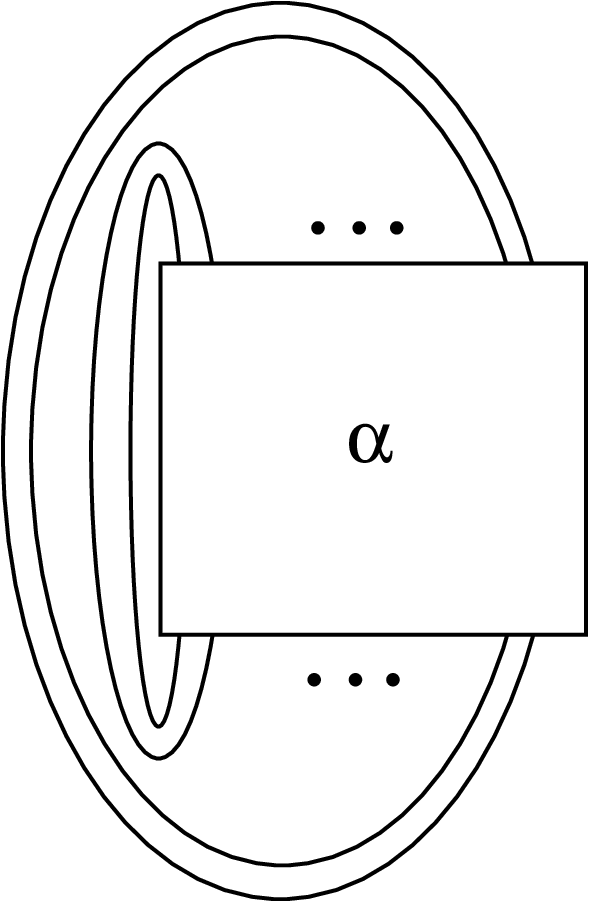}}
\caption{\label{fig:braid-closure1} Closure $\widehat{\a}$ of a braid $\a$}
\end{figure}

It is shown in \cite{Grigorchuk} that the full braid group $B_n$ admits
infinitely many linearly independent homogeneous quasi-morphisms
for every integer $n>2$. However none of these quasi-morphisms are
constructed geometrically. In \cite{surfaces} Gambaudo and Ghys
gave an explicit geometric construction of a family of quasi-morphisms
$\textbf{sign}_n$ on groups $B_n$ defined as follows:
$$\textbf{sign}_n(\a):=sign(\widehat{\a}),$$ where $sign$ is a signature link invariant (see Section \ref{sec-examples}) and $\widehat{\a}$ is the link in $S^3$ which is obtained in the natural way from $\a$ (see Figure \ref{fig:braid-closure1}). In this section we show that knot invariants of certain type define quasi-morphisms on $B_n$ in a similar way.

%%%%%%%%%%%%%%%%%%%%%%%%%%%%%%%%%%%%%%%%%%
%%%%%%%%%%%%%%%%%%%%%%%%%%%%%%%%%%%%%%%%%%

\begin{rem} Another family of quasi-morphisms on $B_n$ (one for each $n$) was constructed recently by Malyutin in \cite{Mal}. These quasi-morphisms are constructed using different methods, in particular they are not defined by knot/link invariants, and will not be discussed in this paper.
\end{rem}

Let us recall some useful notions from knot theory. Let $K$ be a knot in $S^3$. The \textit{four-ball genus} $g_4(K)$ of $K$ is the minimal genus of an oriented surface with boundary which is smoothly embedded in
$\D^4$ such that the image of its boundary under this embedding is the knot $K\subset S^3=\partial\D^4$. We denote by $K^*$ the mirror image of $K$, and by $-K$ the same knot $K$ with the reversed orientation. For any knot $K'$ the knot $K\#K'$
represents the connected sum of $K$ and $K'$. A knot $K$ is called \textit{slice}
if $g_4(K)=0$. We say that $K$ is concordant to $K'$, if there exists a smooth embedding
$\iota:S^1\times[0,1]\hookrightarrow S^3\times[0,1]$ such that
$\iota(S^1\times\{0\})=K$ and $\iota(S^1\times\{1\})=K'$. We
denote by $Conc(S^3)$ the Abelian group whose elements are
equivalence classes of knots in $S^3$, and the multiplication is
the connected sum operation. It is easy to see that the
multiplication is well-defined. The following lemma is a
well-known fact in knot theory, see e.g. \cite{Likorish, Rolfsen}.

\begin{lem}\label{lem:concordance}
For any knot $K$ the knot $K\#-K^*$ is a slice knot.
\end{lem}

\textbf{Main Theorem.} One of the ways to find a quasi-morphism on $B_n$  is to find an
$\R-$valued invariant $I$ of (isotopy classes of) links in $S^3$ so that
$$|I(\widehat{\a\b})-I(\widehat{\a})-I(\widehat{\b})|\leq K_I,$$ where $K_I\geq 0$ depends only on $I$. We will now describe certain ways of closing braids into knots.
\begin{lem}\label{lem:braid-existance}
Let $\b\in B_n$. Then there exists a braid
$\a_{\b}\in B_n$ which satisfies the following properties:\\
\textbf{1.} The closure of $\a_\b\b$ is a knot.\\
\textbf{2.} The closure of $\a_\b$ is a $k$-component unlink for some $1\leq k\leq n$.
\end{lem}
We will say that such a braid $\a_\b$ is a {\it completing braid for $\b$}. Let $\{\s_i\}_{i=1}^{n-1}$ be the standard (Artin) generators of $B_n$.

\begin{proof}
Let $\b\in B_n$. Then there exists a braid $\a$ such that
$\widehat{\a{\b}}$ is a knot and $\widehat{\a}$ is a
$k$-component link for some $1\leq k\leq n$. We write $\a$ as
a product of generators $\s_i^{\pm1}$ and perform a sequence
of crossing changes (i.e. replacing {\it some} of the generators
in the product by their inverses) until we are left with a braid
$\a_\b$ such that $\widehat{\a{\b}}=
\widehat{\a_\b {\b}}$ is a knot and
$\widehat{\a_\b}$ is the $k$-component unlink.
\end{proof}

\begin{defn}\label{defn:alpha-beta}
\rm{Let $I$ be a real-valued knot invariant. Let us fix some
choices of completing braids $\a_\b$ for every $\b\in B_n$. Define a function
$$\widehat{I}:B_n\to\R,$$ by
$\widehat{I}(\b):=I(\widehat{\a_\b\b})$.}
\end{defn}

We will show that under certain conditions $\widehat{I}$ is a quasi-morphism.

\begin{rem}\label{rem:difference-qm} \rm{Note that our function
$\widehat{I}:B_n\to\R$ is slightly different from
an analogous quasi-morphism defined by Gambaudo-Ghys. Their
definition ($\widehat{\textbf{I}}(\b):=I(\widehat{\b})$)
requires $I$ to be a real-valued {\it link} invariant, but in our
definition we require $I$ to be only a real-valued {\it knot}
invariant. For example the Ozsvath-Szabo $\tau$ invariant is
defined only for knots. In case $I$ is a real-valued link
invariant defining a quasi-morphism $\widehat{I}$, the
quasi-morphisms $\widehat{I}$ and $\widehat{\textbf{I}}$ differ by
a constant which depends only on $n$, hence their homogenizations
are equal.}
\end{rem}

\begin{thm} \label{thm:I-quasimorphism}
Suppose that a real-valued knot invariant $I$ defines a
homomorphism $$I:Conc(S^3)\to\R,$$ such that
$|I(K)|\leq c_1g_4(K)$, where  $c_1$ is a real positive constant
independent of $K$. Then $\widehat{I}$ is a quasi-morphism on
$B_n$. Moreover, for a different set of completing
braids $\a_\b$ for every $\b\in B_n$ we get a (possibly
different) quasi-morphism on $B_n$ such that the absolute value of
its difference with $\widehat{I}$ is bounded from above by a
constant depending only on $n$ and therefore the homogenizations
of the two quasi-morphisms are equal.
\end{thm}

\begin{proof}
Take any $\b,\g\in B_n$ and let $\a_\b,\a_\g,
\a_{\b\g}$ be the chosen completing braids for
$\b,\g$ and $\b\g$.

\begin{lem}
\label{lem-cobordism} There exists a cobordism $S$ between the knots
$$(-\widehat{\a_{\b\g}\b\g})^{*}\#(\widehat{\a_\b\b}\#
\widehat{\a_\g\g})\qquad and\qquad(-\widehat{\a_{\b\g}\b\g})^{*}\#\widehat{\a_{\b\g}\b\g}$$
such that $\chi (S)\geq -6n$.
\end{lem}

\begin{proof}
If $T$ is a cobordism between two links $L$ and $L'$, we will
write $L\ \ \overset{T}{\sim}\ \ L'$.

By an observation of Baader (see \cite[Section 4]{Baader}) for any braids $\mu, \nu\in B_n$
$$\widehat{\mu\nu}\ \ \overset{T}{\sim}\ \ \widehat{\mu}\sqcup\widehat{\nu},$$
where $\chi(T)=-n$. Also note that since
$\widehat{\a_{\b\g}},
\widehat{\a_\b},\widehat{\a_\g}$ are unlinks with
no more than $n$ components, we have
$$\widehat{\a_{\b\g}}\ \ \overset{S_1}{\sim}\ \ \widehat{\a_\b}\sqcup\widehat{\a_\g},$$
where $\chi (S_1)\geq 1-2n$. Therefore
$$\widehat{\a_{\b\g}\b\g}\ \ \overset{S_2}{\sim}\ \
\widehat{\a_{\b\g}}\sqcup\widehat{\b}\sqcup\widehat{\g}
\ \ \overset{S_3}{\sim}\ \
\widehat{\a_\b}\sqcup\widehat{\a_\g}\sqcup\widehat{\b}\sqcup\widehat{\g}\
\ \overset{S_4}{\sim}\ \ \widehat{\a_\b\b}\sqcup
\widehat{\a_\g\g}\ \ \overset{S_5}{\sim}\ \
\widehat{\a_\b\b}\# \widehat{\a_\g\g},$$
where $\chi (S_2) = -2n$, $\chi (S_3) \geq 1-2n$ (the cobordism
$S_3$ is the disjoint union of $S_1$ and the trivial cobordism
over $\widehat{\b}\sqcup\widehat{\g}$ given by a disjoint
union of cylinders), $\chi (S_4)=-2n$ and $\chi (S_5)=-1$ ($S_5$ is a saddle cobordism between the disjoint union of two
knots and their connected sum). Thus
$$\widehat{\a_{\b\g}\b\g}\ \ \overset{S_6}{\sim}\ \ \widehat{\a_\b\b}\#
\widehat{\a_\g\g}$$ and hence
$$(-\widehat{\a_{\b\g}\b\g})^{*}\#\widehat{\a_{\b\g}\b\g}\ \
\overset{S_7}{\sim}\ \ (-\widehat{\a_{\b\g}\b\g})^{*}\#(\widehat{\a_\b\b}\#
\widehat{\a_\g \g}),$$ where
$$\chi (S_6) = \chi(S_7) = \chi(S_2) + \chi(S_3) +\chi (S_4) +\chi
(S_5) \geq -6n,$$ as required.
\end{proof}

Let us now finish the proof of the theorem. Lemma
\ref{lem:concordance} implies that the knot
$$(-\widehat{\a_{\b\g}\b\g}^*)\#\widehat{\a_{\b\g}\b\g}$$ is
a slice knot. Therefore using Lemma \ref{lem-cobordism} we get
$$g_4(-(\widehat{\a_{\b\g}\b\g}^*)\#(\widehat{\a_\b\b}\#\widehat{\a_\g\g}))\leq 3n.$$
This yields
$$|I(-(\widehat{\a_{\b\g}\b\g}^*)\#(\widehat{\a_\b\b}\#\widehat{\a_\g\g}))|
\leq c_1\cdot3n.$$ Applying  the equalities
\begin{equation}\label{eqn-properties-of-I}
I(-K^*)=-I(K),\quad\quad\quad I(K\#K')=I(K)+I(K'),
\end{equation}
we get
$$|\widehat{I} (\b\g) - \widehat{I}(\b) -
\widehat{I}(\g)|\leq 3c_1n,$$ which means that $\widehat{I}$ is
a quasi-morphism.

Finally note that for any two choices $\a_\b$ and
$\a'_\b$ of completing braids for $\b\in B_n$ one can
show, similarly to the proof of Lemma \ref{lem-cobordism}, that
$$g_4 ((\widehat{\a_\b
\b})\# -(\widehat{\a'_\b \b})^*)\leq c_2 n,$$ for
some positive constant $c_2$ independent of $\b, \a_\b,
\a'_\b$. Hence by equations ~\eqref{eqn-properties-of-I} we get
$$|I(\widehat{\a_\b \b}) -
I(\widehat{\a'_\b \b})| = |I((\widehat{\a_\b
\b})\# -(\widehat{\a'_\b \b})^*)| \leq c_1\cdot c_2n.$$
Thus different choices of completing braids yield (possibly
different) quasi-morphisms on $B_n$ whose difference is bounded in
absolute value from above by a constant depending only on $n$ and
therefore their homogenizations are equal.
\end{proof}

\begin{rem}
\rm{Note that $g_4$ defines a semi-norm on $Conc(S^3)$ and Theorem \ref{thm:I-quasimorphism} can be reformulated as follows: each element of $\Hom(Conc(S^3),\R)$, which is Lipshitz with respect to the semi-norm defines a quasi-morphism on $B_n$.}
\end{rem}

\section{Examples, properties and applications}

\label{sec-examples}
In this section we discuss the following knot/link invariants: link $\o$-signatures, the Rasmussen link invariant $s$ and the Ozsvath-Szabo knot invariant $\tau$.

First, let us recall a few definitions and notations concerning
braids and quasi-morphisms. The \textit{braid length} $l(\g)$ of $\g\in B_n$ is the length of the shortest word representing $\g$ with respect to the generators $\s_{1},\ldots,\s_{n-1}$. Throughout the paper the induced homogeneous quasi-morphism obtained by the homogenization of a quasi-morphism
$\varphi$ (see ~\eqref{eq:definition-homogenization}) will be denoted by $\widetilde{\varphi}$.

\subsection{Signature and $\o$-signature quasi-morphisms}
\label{subsec-signatures}
Let $L$ be an oriented link in $S^3$, then
there exists an oriented surface $\Sigma_L$ with boundary $L$. It is called a \textit{Seifert surface} of $L$.
We choose a basis $\{b_1,\ldots,b_{2g+|L|-1}\}$ in
$H_1(\Sigma_L,\Z)$ and define a symmetric bilinear form on
$H_1(\Sigma_L,\Z)$ as follows:
$$\Omega(b_i, b_j)=\lk(b_i,b_j^+)+\lk(b_j,b_i^+),$$
where $\lk$ is the linking number and $b_i^+$ is a push-off of the
curve, which represents $b_i$ in $\Sigma_L$, from $\Sigma_L$ along the
positive normal direction  to $\Sigma_L$.
Tensoring by $\R$ we get a symmetric bilinear form on
$H_1(\Sigma_L,\R)$. The signature of this form is
independent of the choices of $\Sigma_L$ and the
basis of $H_1(\Sigma_L,\Z)$, see \cite{Likorish,M-sign}.
Thus it is an invariant of $L$ and is denoted by $sign(L)$.

For any complex number $\o\neq 1$ and link $L$ there exists the
$\o$-signature link invariant $sign_\o(L)$, such that
$sign_{-1}(L)=sign(L)$. It is defined as follows. We tensor the
bilinear form $\Omega$ by  $\mathbb{C}$, and we get a bilinear
form on $H_1(\Sigma_L,\mathbb{C}).$
The signature of the following hermitian form
$$\Omega_{\o}(b_i,b_j)=(1-\o)\Omega(b_i,\overline{b_j})+
(1-\overline{\o})\Omega(\overline{b_j},b_i)$$
on $H_1(\Sigma_L,\mathbb{C})$ is independent of the choice of the Seifert surface $\Sigma_L$ and the
basis for $H_1(\Sigma_L,\Z)$ (see e.g. \cite{Likorish}).
Thus $sign_\o(L)$ is a link invariant for each $\o\neq 1$.

In \cite{w-signature} Gambaudo and Ghys showed that $sign_\o$ defines a quasi-morphism on $B_n$. Alternatively, $sign_\o$ satisfies conditions of Theorem \ref{thm:I-quasimorphism} for each $\o$ (see e.g. \cite{Likorish}). Hence
\begin{cor}
Both $\widehat{sign}:B_n\to\R$ and $\widehat{sign}_\o:B_n\to\R$ are quasi-morphisms.
\end{cor}

The induced homogeneous quasi-morphisms on $B_n$ are denoted by $\widetilde{sign}$ and $\widetilde{sign}_\o$ respectively. Following \cite{surfaces} we denote by $lk:B_n\to\Z$ the unique (up to the multiplication by a constant) homomorphism from $B_n$ to $\Z$ by setting $lk(\s_i^{\pm1})=\pm 1$. Let us recall the following

\begin{defn}\label{defn:alternating-link}\rm
A link diagram is called \textit{alternating} if the crossings alternate under, over, under, over, and so on as one travels along each component of the link. A link is called \textit{alternating} if it has an alternating diagram. A braid $\a\in B_n$ is called \textit{alternating} if its closure $\widehat{\a}$ is an alternating link diagram.
\end{defn}

\begin{prop}\label{prop:alt-braid}
Let $\g\in B_n$ be an alternating braid, then
$$\widetilde{sign}(\g)=lk(\g).$$
\end{prop}
The proof of this proposition will be given in the next subsection.
\begin{rem}\rm
One can easily show that $\widetilde{sign}$ is not a homomorphism for each $n>2$.
\end{rem}

\subsection{Rasmussen quasi-morphism}

The Rasmussen invariant $s(K)$ of a knot $K$ in $S^3$ was
discovered by Jacob Rasmussen in 2004, see \cite{Rasmussen}. It comes from the Lee theory \cite{Lee} which is closely related to the Khovanov homology \cite{Khovanov}. This is a very powerful knot invariant which was used by Rasmussen in \cite{Rasmussen} to give a first combinatorial proof of the Milnor conjecture. This invariant was extended to links in \cite{Beliakova}.

\begin{figure}[htb]
\centerline{\includegraphics[height=0.6in]{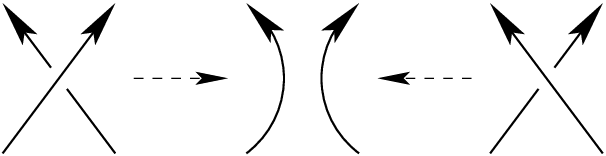}}
\caption{\label{fig:Conway-smoothing} Smoothing of a crossing}
\end{figure}

Before we list the properties of $s$ let us recall a notion of the Seifert algorithm.
Let $D_L$ be a diagram of an oriented link $L$.
Let us smoothen each crossing in $D_L$ as shown in Figure
\ref{fig:Conway-smoothing}.
%%%%%%%%%%%%%%%%%%%%%%%%%%%%%%%%%%%%%%%%%
%%%%%%%%%%%%%%%%%%%%%%%%%%%%%%%%%%%%%%%%%
 The resulting smoothed diagram $\widehat{D}_L$
consists of oriented simple closed curves, which are called
the \textit{Seifert circles}. Thus $\widehat{D}_L$ is the boundary
of a union of disjoint discs. We join these discs together with
half-twisted strips corresponding to the crossings in the diagram.
This yields an oriented surface bounded by $L$.
If this surface is disconnected, then we connect its components by the connected sum operation.
The above algorithm is called the \textit{Seifert algorithm}.

We denote by $K_+$ and $K_-$ the knots which differ by a single crossing change:
from a positive crossing in $K_+$ to a negative one in $K_-$ (positive/negative
crossings are shown on the left/right in Figure \ref{fig:Conway-smoothing}). Now we list eight properties of the Rasmussen link invariant $s$. The first seven were proved in \cite{Rasmussen} and the last one was proved in \cite{Baader}.
\begin{equation}\label{eq:property1}
s(K)=s(-K),
\end{equation}
\begin{equation}\label{eq:property2}
s(K\# K')=s(K)+s(K'),
\end{equation}
\begin{equation}\label{eq:property-mirror-of-K}
s(K^*)=-s(K),
\end{equation}
\begin{equation}\label{eq:property3}
|s(K)|\leq 2g_4(K),
\end{equation}
\begin{equation}\label{eq:property4}
s(K_-)\leq s(K_+)\leq s(K_-)+2.
\end{equation}
For an alternating knot $K$ we have
\begin{equation}\label{eq:property7}
s(K)=sign (K).
\end{equation}
For a positive knot $K$ we have
\begin{equation}\label{eq:property5}
s(K)=w(D_K)-o(D_K)+1,
\end{equation}
where $D_K$ is some positive diagram of $K$, $w(D_K)$ is the
number of positive crossings in $D_K$ minus the number of
negative crossings in $D_K$, and $o(D_K)$ is the number of
Seifert circles. For any knot $K$ we have
\begin{equation}\label{eq:property6}
1+w(D_K)-o(D_K) \leq s(K)\leq -1+w(D_K)+o(D_K),
\end{equation}
where $D_K$ is any diagram of $K$. It is shown in
\cite{Baader} that ~\eqref{eq:property6} follows from ~\eqref{eq:property-mirror-of-K},
~\eqref{eq:property4} and ~\eqref{eq:property5}.

In \cite{Baader} Baader proved that the Rasmussen link invariant $s$ induces a quasi-morphism  on the
braid group $B_n$. Alternatively, note that properties ~\eqref{eq:property1}, ~\eqref{eq:property2}, ~\eqref{eq:property-mirror-of-K} and ~\eqref{eq:property3} imply that $s$ satisfies conditions of Theorem \ref{thm:I-quasimorphism} and hence
$\widehat{s}:B_n\to\R$ is a quasi-morphism.

\begin{thm}\label{thm:link-ras}
Suppose that $\widehat{I}: B_n\to\R$ is a
quasi-morphism defined by a real-valued knot invariant $I$,
which satisfies property \eqref{eq:property6}, where $s$ is
substituted by $I$. Then $\widetilde{I}=lk$.
\end{thm}
\begin{proof}
Take $\b\in B_n$. Then for all $p>0$ there exists a completing
braid $\a_p$ in $B_n$, such that
$\widehat{\a_p\b^p}$ is a knot and $|lk(\a_p)|\leq M(n)$, where $M(n)$ is some real positive constant which depends only on $n$.
For example if $\b\in P_n$, than we can take
$\a_p=\s_1\cdot\ldots\cdot\s_{n-1}$ for all $p\in\N$. Note that
$$w(D_{\widehat{\a_p\b^p}})=lk(\a_p\b^p)\quad \textrm{and}\quad o(D_{\widehat{\a_p\b^p}})=n.$$
The second equation follows from the fact that for any chosen orientation on the knot $\widehat{\a_p\b^p}$ \textbf{all} strands in $\a_p\b^p$ are oriented from the bottom to the top, or from the top to the bottom (this depends on the orientation) and hence the number of Seifert circles is $n$. Property \eqref{eq:property6} yields
\begin{equation*}
1+lk(\a_p\b^p)-n\leq I(\widehat{\a_p\b^p})\leq -1+lk(\a_p\b^p)+n.
\end{equation*}
Hence
\begin{equation*}
\lim\limits_{p\to\infty}\frac{1}{p}(1+lk(\a_p\b^p)-n)
\leq \lim\limits_{p\to\infty}\frac{1}{p}I(\widehat{\a_p\b^p})
\leq \lim\limits_{p\to\infty}\frac{1}{p}(-1+lk(\a_p\b^p)+n).
\end{equation*}
Therefore
$\widetilde{I}(\b)=\lim\limits_{p\to\infty}\frac{1}{p}\widehat{I}(\b^p)=\lim\limits_{p\to\infty}\frac{1}{p}I(\widehat{\a_p\b^p})=lk(\b)$.
\end{proof}

\begin{cor}\label{cor:Ras-lk}
Each $\a\in B_n$ satisfies $\widetilde{s}(\a)=lk(\a)$.
\end{cor}

Now we are ready to prove Proposition \ref{prop:alt-braid}.

\begin{proof}[Proof of Proposition \ref{prop:alt-braid}]
For each $p\in\N$ the braid $\g^p$ is alternating and there exists an
alternating braid $\a_p$ such that $l(\a_p)<M(n)$, where $l$ is the braid length, $\widehat{\a_p\g^p}$ is an alternating knot and $M(n)$ is some real positive constant which depends only on $n$. Hence
$$\widetilde{sign}(\g)=\lim\limits_{p\to\infty}\frac{sign(\widehat{\a_p\g^p})}{p}=
\lim\limits_{p\to\infty}\frac{s(\widehat{\a_p\g^p})}{p}=\lim\limits_{p\to\infty}\frac{s(\widehat{\g^p})}{p}=\widetilde{s}(\g)=lk(\g),$$
where the second equality follows from \eqref{eq:property7}.
\end{proof}

In \cite{Baader} Baader  defined the following
quasi-morphism:
$$s-lk+n-1:B_n\to\R,$$
where $s(\b):=s(\widehat{\b})$ (Gambaudo-Ghys definition).
This quasi-morphism maps all positive braids to zero. Therefore it
descends to a quasi-morphism on $B_n/\langle\Delta_n\rangle$,
where $\Delta_n=(\s_1\s_2\cdot\ldots\cdot\s_{n-1})^n$. In \cite{Baader} Baader asked whether this quasi-morphism is bounded for $n>2$. An immediate corollary of Theorem \ref{thm:link-ras} gives a positive answer to this question.

\begin{cor}
The quasi-morphism $s-lk+n-1:B_n\to\R$ is bounded.
\end{cor}

\subsection{Connection between $\tau$, $s$ and the braid index}

In \cite{Ozsvath-tau} Ozsvath and Szabo
defined a knot invariant $\tau$ which comes from the knot Floer homology \cite{MOS, Ozsvath, Rasmussen-knot-floer}.
This invariant satisfies in particular the following properties:
\begin{equation}\label{eq:property-tau1}
\tau(K^*)=-\tau(K),
\end{equation}
\begin{equation}\label{eq:property-tau2}
\tau(K\#K')=\tau(K)+\tau(K'),
\end{equation}
\begin{equation}\label{eq:orientation-tau}
\tau(-K)=\tau(K),
\end{equation}
\begin{equation}\label{eq:four-ball-for-tau}
|\tau(K)|\leq g_4(K),
\end{equation}
\begin{equation}\label{eq:property-tau3}
0\leq\tau(K_+)-\tau(K_-)\leq 1,
\end{equation}
\begin{equation}\label{eq:property-tau4}
\tau(\widehat{\a})=\frac{lk(\a)-n+1}{2},
\end{equation}
where $\a\in B_n$ is a positive braid (i.e., it may be written as a product of positive powers of standard generators $\s_1,\ldots,\s_{n-1}$ of $B_n$). The first five properties were proved in \cite{Ozsvath-tau}, and the last one was
proved in \cite{Livingston}. The invariants $s$ and $2\tau$ share similar properties and coincide on positive and alternating knots. It was conjectured by Rasmussen  \cite{Rasmussen} that they are equal. This conjecture was disproved by Hedden and Ording \cite{HO}. In this subsection we show that $\widetilde{s}=2\widetilde{\tau}$, and that the absolute value of the difference between $2\tau$ and $s$ is bounded by twice the braid index of the knot.

\begin{lem}\label{lem:tau-property}
Let $\b\in B_n$ such that $\widehat{\b}$ is a knot. Then
\begin{equation}\label{eq:property-tau}
lk(\b)-n+1 \leq 2\tau(\widehat{\b})\leq lk(\b)+n-1.
\end{equation}
\end{lem}
\begin{proof}
The proof of the lower bound in \eqref{eq:property-tau} is exactly
the same as in \cite{Shumakovitch}, where $s$
is substituted by $2\tau$. We present it for the reader's
convenience. Let $K$ be any knot in $S^3$. Recall that
\begin{equation}\label{eq:tau+-}
0\leq 2\tau(K_+)-2\tau(K_-)\leq 2,
\end{equation}
\begin{equation}\label{eq:tau-mirror}
\tau(K^*)=-\tau(K).
\end{equation}
Let $\a\in B_n$ be any positive braid such that $\widehat{\a}$ is a knot. It follows from property ~\eqref{eq:property-tau4} that the following equality holds:
\begin{equation}\label{eq:tau-positive}
2\tau(\widehat{\a})=lk(\a)-n+1.
\end{equation}
It means that the left inequality in \eqref{eq:property-tau} is
then an equality. When a positive crossing of $\widehat{\b}$ is changed
into a negative one, the number $lk(\b)-n+1$ decreases by
$2$, while $2\tau(\widehat{\b})$ decreases by at most $2$, because of
\eqref{eq:tau+-}. Hence, the left inequality in
\eqref{eq:property-tau} is preserved. The upper bound in
\eqref{eq:property-tau} is true, because of \eqref{eq:tau-mirror}
combined with the lower bound.
\end{proof}

\begin{cor}
The knot invariant $\tau$ defines a quasi-morphism on $ B_n$ and $$2\widetilde{\tau}=lk.$$
\end{cor}

\begin{proof}
Note that properties ~\eqref{eq:property-tau1}, ~\eqref{eq:property-tau2}, ~\eqref{eq:orientation-tau} and ~\eqref{eq:four-ball-for-tau} imply that $\tau$ satisfy conditions of Theorem \ref{thm:I-quasimorphism}, and hence $\widehat{\tau}:B_n\to\R$ is a quasi-morphism.
The second statement follows from Lemma \ref{lem:tau-property} and the proof of Theorem \ref{thm:link-ras}.
\end{proof}

Recall that every knot $K$ in $\R^3$ can be presented as a closure of some braid in $B_n$. The braid index of $K$ is the minimal such $n$. It is denoted by $br(K)$.

\begin{thm}
Let $\widehat{s}$ and $\widehat{\tau}$ be the quasi-morphisms on $B_n$ which are induced from Rasmussen and Ozsvath-Szabo knot invariants. Then for every $\b\in B_n$ the following inequality holds
$$|\widehat{s}(\b)-2\widehat{\tau}(\b)|\leq 2(n-1).$$
\end{thm}

\begin{proof}
As we explained before both $s$ and $2\tau$ satisfy property
\eqref{eq:property-tau}. It means that
$$|\widehat{s}(\b)-lk(\a_\b\b)|\leq n-1\quad \textrm{and}\quad |2\widehat{\tau}(\b)-lk(\a_\b\b)|\leq n-1.$$
Thus by the triangle inequality we have
$$|\widehat{s}(\b)-2\widehat{\tau}(\b)|\leq 2(n-1).$$
\end{proof}

\begin{cor}
For every knot $K$ the following inequality holds: $$|s(K)-2\tau(K)|\leq 2(br(K)-1).$$
\end{cor}

\begin{pr*}
Does there exist a real-valued knot invariant $I$,
independent from the $\o$-signatures, $\tau$ and $s$, such that $I$ satisfies the conditions of Theorem \ref{thm:I-quasimorphism}, and such that the induced homogeneous quasi-morphism $\widetilde{I}:B_n\to\R$ is non-trivial?
\end{pr*}

\section{Induced quasi-morphisms on $\cD$}

\subsection{Gambaudo-Ghys construction}
Recall that $\cD:=Diff^{\infty}(\D^{2}, \partial\D^{2}, \area)$ is the group
of compactly supported area-preserving diffeomorphisms of the open unit disc $\D^2$. In this subsection we discuss the construction, essentially contained in \cite{surfaces}, which takes a homogeneous quasi-morphism on $P_n$ and produces from it a quasi-morphism on $\cD$.

Denote by $X_n$ the space of all \emph{ordered} $n$-tuples of distinct points in $\D^2$. Let us fix a base point $\overline{z}=(z_1,\ldots,z_n)\in X_n$ and let $\overline{x}=(x_1,\ldots, x_n)$ be any other point in $X_n$. It is a well-known fact that $\cD$ is path-connected (see e.g. \cite{McDuff}).
Take  $g\in\cD$ and  any path $\psi_t$, $0 \leq t\leq 1$, in $\cD$
between $Id$ and $g$. Connect $\overline{z}$ to $\overline{x}$ by a straight line in $(\D^2)^n$, then act on $\overline{x}$ with the path $\psi_t$, and then connect $g(\overline{x})$ to $\overline{z}$ by the straight line in $(\D^2)^n$. We get a loop in $(\D^2)^n$. More specifically it looks as follows. Connect $z_i$ to $x_i$ by straight lines $\mathfrak{l}_{1,i}:\left[0,\frac{1}{3}\right]\to\D^2$
in the disc, then act with the path $\psi_{3t-1}$, $\frac{1}{3} \leq t\leq \frac{2}{3}$, on each $x_i$, and then
connect  $g(x_i)$ to $z_i$ by straight lines $\mathfrak{l}_{2,i}:\left[\frac{2}{3},1\right]\to\D^2$
in the disc, for all $1\leq i \leq n$. Note that for almost all $n$-tuples of different points $x_{1},\ldots, x_{n}$ in the disc the concatenations of the paths $\mathfrak{l}_{1,i}:\left[0,\frac{1}{3}\right]\to\D^{2}$, $\psi_{3t-1}:\left[\frac{1}{3},\frac{2}{3}\right]\to\D^{2}$ and $\mathfrak{l}_{2,i}:\left[\frac{2}{3},1\right]\to\D^{2}$, $i=1,\ldots,n,$ yield a loop in $X_{n}$. The homotopy type of this loop is an element in $P_n$ (here $P_n$ is identified with the fundamental group $\pi_1 (X_n,\overline{z})$). This element is independent of the choice of $\psi_t$ because $\cD$ is contractible (see e.g. \cite{contractible-group,T}), it will be denoted by $\g(g;\overline{x})$. Let $\widetilde{\varphi}_n$ be a homogeneous quasi-morphism on $P_n$. Denote $d\overline{x}:=dx_1\cdot\ldots\cdot dx_n$ and set
\begin{equation}
\Phi (g)=\int\limits_{X_n}\widetilde{\varphi}_n(\g(g;\overline{x}))d\overline{x}.
\end{equation}

\begin{lem}[c.f. \cite{surfaces}] \label{lem:integrability}
The function $\Phi$ is well defined. It is
a  quasi-morphism on $\cD$.
\end{lem}

As an immediate corollary we get that the formula
$$\widetilde{\Phi} (g)=\lim\limits_{p\to\infty}\frac{1}{p}
\int\limits_{X_n}\widetilde{\varphi}_n(\g(g^p;\overline{x}))d\overline{x}$$ yields a  homogeneous
quasi-morphism $\widetilde{\Phi}: \cD\to \R$.

\begin{proof}

Let $g\in\cD$. For any isotopy $\psi_t$, $0\leq t\leq 1$, in $\cD$ between $Id$
and $g$, any $\overline{x}\in X_n$ and any $1\leq i,j\leq n,\thinspace i\neq j$ denote
$$l_{i,j}(t):=\frac{\psi_t(x_i)-\psi_t(x_j)}{\|\psi_t(x_i)-\psi_t(x_j)\|}:[0,1]\to S^1\quad \textrm{and} \quad  L_{i,j}(\overline{x}):=\frac{1}{2\pi} \int\limits_{0}^{1}\left\|\frac{\partial}{\partial t}(l_{i,j}(t))\right\|dt,$$
where $\|\cdot\|$ is the Euclidean norm. Note that $L_{i,j}(\overline{x})$
is the length of the path $l_{i,j}(t)$ divided by $2\pi$. It follows that $L_{i,j}(\overline{x})+4$ is an upper bound for the number of times the string $i$ turns around the string $j$ in the positive direction plus the number of times the string $i$ turns around the string $j$ in the negative direction in the braid $\g(g;\overline{x})$. Thus we get the following inequality
\begin{equation}\label{eq:inequality1}
\displaystyle\sum_{i<j}^{n}2\left(L_{i,j}(\overline{x})+4\right)\geq l(\g(g;\overline{x})).
\end{equation}
Take any finite generating set $\cS$ of $P_n$. Note that for any homogeneous quasi-morphism
$\widetilde{\varphi}_n:P_n\to\R$ one has
\begin{equation}\label{eq:length-inequality}
|\widetilde{\varphi}_n(\g)|\leq\left(D_{\widetilde{\varphi}_n}+
\max\limits_{\xi\in \cS}|\widetilde{\varphi}_n(\xi)|\right)l_\cS(\g),
\end{equation}
where $l_\cS(\g)$ is the length of a word $\g$ with respect to $\cS$. Recall that $l(\g)$ is the length of $\g$ with respect to the set $\{\s_i\}_{i=1}^{n-1}$. It follows from \cite[Corollary 24]{Harp} that there exist two positive constants $K_{1,\cS}$ and $K_{2,\cS}$, which are independent of $\g$,
such that $$l_\cS(\g)\leq K_{1,\cS}\cdot l(\g)+K_{2,\cS}.$$

It follows from \eqref{eq:length-inequality} that
\begin{equation}\label{eq:inequality2}
|\widetilde{\varphi}_n(\g(g;\overline{x}))|\leq N_1l(\g(g;\overline{x}))+N_2,
\end{equation}
where $N_1 =K_{1,\cS}(D_{\widetilde{\varphi}_n}+\max\limits_{\xi\in\cS}|\widetilde{\varphi}_n(\xi)|)$ and
$N_2=K_{2,\cS}(D_{\widetilde{\varphi}_n}+\max\limits_{\xi\in\cS}|\widetilde{\varphi}_n(\xi)|)$. Inequalities
\eqref{eq:inequality1} and \eqref{eq:inequality2} yield the
following inequality:
\begin{equation*}
|\widetilde{\varphi}_n(\g(g;\overline{x}))|\leq 2N_1\left(\displaystyle\sum_{i<j}^{n}L_{i,j}(\overline{x})+4\right)+N_2.
\end{equation*}
It follows that
\begin{equation*}
\left|\Phi(g)\right|\leq 2N_1 \cdot
\vol((\D^2)^{n-2})\left(\sum\limits_{i<j}^n
\int\limits_{\D^2\times\D^2}L_{i,j}(\overline{x}) dx_idx_j\right)+(4N_1n(n+1)+N_2)\vol((\D^2)^n).
\end{equation*}
It follows from \cite[Lemma 1]{Gambaudo-Lagrange} that the integral in the above inequality
is well defined and hence $\left|\Phi(g)\right|<\infty$. The proof of the fact that $\Phi$ is a quasi-morphism is exactly the same as in \cite{surfaces}.
\end{proof}

From now on the homogeneous quasi-morphism on $\cD$ induced by $\widetilde{\varphi}_n:P_n\to\R$ will be denoted by $\widetilde{\Phi}$. We denote by $\C:\cD\to\R$ the celebrated Calabi homomorphism (see \cite{Banyaga, Calabi}, cf. \cite{enlacement}).

\begin{rem}\rm  It follows from the interpretation of $\C$ in \cite{enlacement} that the homogeneous quasi-morphisms on $\cD$ induced by the homogeneous Rasmussen and Ozsvath-Szabo quasi-morphisms $\widetilde{s}$ and $\widetilde{\tau}$ are equal to $\C$ multiplied by $2n(n-1)\pi^{n-1}$ and by\\ $n(n-1)\pi^{n-1}$ respectively.
\end{rem}

\subsection{Generic autonomous Hamiltonians}
\label{subsec-aut-ham}
Denote the space of autonomous compactly-supported Hamiltonians $H:\D^2\to\R$ by $\cH$. In this subsection we define the notion of a Morse-type Hamiltonian. It follows from \cite[Theorem 2.7]{Milnor} that Morse-type
Hamiltonians form a $C^1$-dense subset of $\cH$. For any Hamiltonian $H$ in  this subset we  present a calculation of
$\widetilde{\Phi}$ on the time-one flow of $H$ (Theorem \ref{thm:general-formula}). Our result is presented in terms of
the integral of the push-forward of $H$ to its Reeb graph against a certain signed measure on the graph. Similar formulas on other surfaces were established in \cite{EP,P1,P2}. This result enables us to relate the asymptotic behavior of an induced signature quasi-morphism with the Calabi homomorphism (Theorem \ref{thm:asymptotics}).

\begin{defn}\label{defn:Morse-type-functions}
\rm{We say that a  function $H\in\cH$ is
of  \textit{Morse-type}  if:\\
 \textbf{1.} There exists a connected open
neighborhood $U$ of $\partial\D^2$, such that
$\partial\overline{U}\setminus\partial\D^2$ is a smooth
simple curve,
$H|_{\overline{U}}\equiv 0$ and $H$ has no degenerate critical points in  $\D^2\setminus \overline{U}$.\\
\textbf{2.} There exists an open set $V\supset\overline{U}$ such that $H|_{V\setminus\overline{U}}$ has no critical points.\\
\textbf{3.} The inequality $H(x)\neq H(y)$ holds for each two non-degenerate different critical points $x$ and $y$.}
\end{defn}

\begin{defn}
\rm{\textit{A charged tree} $T$ is a finite tree equipped with a signed measure, such that each edge has a total finite measure.}
\end{defn}

\begin{figure}[htb]
\centerline{\includegraphics[height=1.3in]{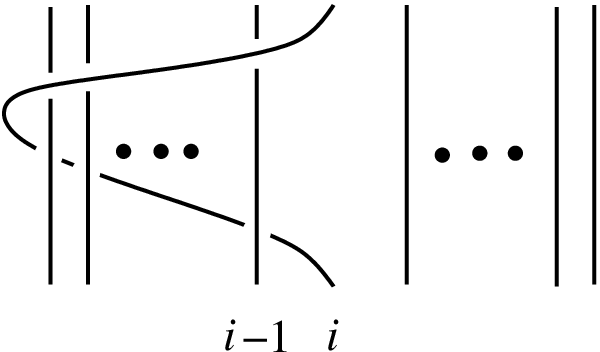}}
\caption{\label{fig:braids-eta-i-n} Braid $\eta_{i,n}$}
\end{figure}

Let $\widetilde{\varphi}_n:P_n\to\R$ be any homogeneous quasi-morphism. Here we explain how we associate to
$(H,\widetilde{\varphi})$ a charged tree $(T,\mu)$.
Let $H\in\cH$ be a Morse-type function with $l$ critical
points in $\D^2\setminus\overline{U}$, where $U$ is as in Definition \ref{defn:Morse-type-functions}, and let $c=\partial\overline{U}\setminus\partial\D^2$.\\
\textbf{Step 1.} Let us recall the notion of angle-action
symplectic coordinates. We remove from $\D^2$ all singular level curves inside $\D^2\setminus\overline{U}$, the curve $c$  and $\partial\D^2$.
We get $l$ open annuli $A_i$ (a punctured disc is also viewed as an annulus) and an open annulus $A_{l+1}=U^\circ$. Each $\D^2\setminus A_i$ has two connected components. We denote by $CA_i$ the component which does not contain $\partial\D^2$ and by
$$a_i:=\partial(CA_i)\cap\partial\overline{A}_i\hspace{10mm} R_i:=\frac{\textrm{area}(A_i)}{2\pi}\hspace{10mm} CR_i:=\frac{\textrm{area}(CA_i)}{2\pi}.$$
By the Liouville-Arnold theorem (see \cite{Arnold}) on each one  of the annuli $A_i$ ($1\leq i\leq l$),
there exist so-called angle-action symplectic coordinates $(\theta_i,J_i)$,
$\theta_i\in [0,2\pi]$, and a $C^\infty$-function
$$\hbar_i:[CR_i,CR_i+R_i]\to\R,$$
such that on each level curve $c_i$ in $A_i$  one has
$H|_{c_i}=\hbar_i(J_i)$, where the coordinate $J_i$ along $c_i$ is equal to the sum of $CR_i$ and
the area, divided by $2\pi$, of the annulus bounded by $c_i$ and $a_i$. Let $h_1$ be the time-one Hamiltonian flow generated by $H$. In these coordinates $h_1$ moves points on each level curve with a constant speed
$$\hbar_i^{\prime}=\frac{\partial \hbar_i}{\partial J_i}.$$
Let $p,q$ be the coordinates on $\R^2$ and note that  $H|_{A_{l+1}}=0$. Then
$$\int\limits_{\D^2}H(p,q)dpdq=\sum\limits_{i=1}^l\int\limits_{A_i}H(p,q)dpdq
=2\pi\sum\limits_{i=1}^l\int\limits_{CR_i}^{CR_i+R_i}\hbar_i(J_i)dJ_i.$$
\textbf{Step 2.} Let $H$ be a
Morse-type function with $l$ critical
points in $\D^2\setminus\overline{U}$. Define the equivalence relation $\sim$ on
points on $\D^2$ by $x\sim y$ whenever $x$, $y$ are in the
same connected component of a level set of $H$. The Reeb graph $T$
corresponding to $H$ is the quotient of $\D^2$ by the
relation $\sim$. The edges of $T$ come from the annuli in
$\D^2$ fibered by level loops, the valency-one vertices
correspond to the min/max critical points of $H$ and to
$\partial\D^2$, and the valency-three vertices correspond
to the saddle critical points of $H$. In our case $T$ is a rooted
tree with the root being the vertex of $T$ corresponding to the
domain $\overline{U}\supset \partial\D^2$. Each edge $e_j$
in $T$ corresponds to an annulus $A_j$ for each $0\leq j\leq l$. A Reeb graph is called \textit{simple} if it has only 2 vertices and one edge.\\
\textbf{Step 3.} The action coordinate on any of the annuli
induces a coordinate $J_j$ on the corresponding edge $e_j$ of $T$.
Thus each Morse-type Hamiltonian $H$ descends to a function
$\hbar:T\to\R$, where $\hbar=\hbar_j$ and
$\hbar^\prime=\hbar_j^\prime$ on each open edge $e_j$. Here
$j\in\{0,\ldots,l\}$. For a homogeneous quasi-morphism
$\widetilde{\varphi}_n:P_n\to\R$ and $0\leq j\leq l$ we define a signed measure
$\mu$ on $T$ by setting
$$d\mu(J_j):=(2\pi)^{n}\displaystyle\sum_{i=2}^n\widetilde{\varphi}_n(\eta_{i,n})i\dbinom{n}{i}(J_j)^{i-1}
\left(\frac{1}{2}-J_j\right)^{n-i} dJ_j,$$
where $J_j$ is the coordinate on each edge $e_j$ and $\eta_{i,n}$ is the pure braid in $P_n$ shown in Figure \ref{fig:braids-eta-i-n}. Note that all of these braids commute with each other.
%%%%%%%%%%%%%%%%%%%%%%%%%%%%%%%%%%%%%%%%%%%%%%%%%%%%%%%%%%%%%%%%%%%%%%%
%%%%%%%%%%%%%%%%%%%%%%%%%%%%%%%%%%%%%%%%%%%%%%%%%%%%%%%%%%%%%%%%%%%%%

\begin{thm}\label{thm:general-formula}
Let $H\in\cH$ be a Morse-type function,
$\widetilde{\varphi}_n:B_n\to\R$ a homogeneous
quasi-morphism, and $\widetilde{\Phi}$ the corresponding homogeneous quasi-morphism on $\cD$. Then
\begin{equation}\label{eq:general-formula}
\widetilde{\Phi}(h_1)=\int\limits_{T}\hbar^\prime d\mu,
\end{equation}
where $h_1$ is the time-one Hamiltonian flow generated by $H$.
\end{thm}

\begin{proof}[Sketch of the Proof]
Let $l$ be the number of isolated critical points of $H$. If $l=1$ then the Reeb graph $T$ which corresponds to $H$ is simple, and Theorem \ref{thm:general-formula} is just \cite[Lemma 5.2]{surfaces}. Let $l>1$. Take a sequence of nonnegative integers $n_1,...,n_{l+1}$ such that $n=n_1+...+n_{l+1}$. Let $\overline{x}\in X_n$ such that $\overline{x}=(x_1,...,x_{n_1},x_{n_1+1},...,x_{n_1+n_2},...,x_{n_1+n_2+...+n_l},...,x_n)$ where $x_i$ is in $A_j$ if $n_1+n_2+...+n_{j-1}+1\leq i\leq n_1+n_2+...+n_j$ for each $1\leq i\leq n$ and $1\leq j\leq l+1$. We also require that for each $1\leq j\leq l$ and $x_i,x_k\in A_j$ such that $i<k$ we have $H(x_i)<H(x_k)$. Then for every $p\in\N$
$$\g(h_1^p,\overline{x})=\a_{1,p}\g_{1,p}\cdot\ldots\cdot\g_{l+1,p}\a_{2,p},$$
where $|\widetilde{\varphi}_n(\a_{1,p})|$ and $|\widetilde{\varphi}_n(\a_{2,p})|$ are bounded by a constant independent of $\overline{x}$ and $p$, each braid $\g_{j,p}$ depends up to conjugation only on $p$ and points $x_{n_1+n_2+...+n_j+1},...,x_{n_1+n_2+...+n_j+n_{j+1}}$. All of $\g_{j,p}$ commute with each other and up to conjugation are elements of the Abelian group which is generated by the braids $\eta_{i,n}$. Now we apply angle-action coordinates in each annuli $A_j$ and proceed as in the proof of \cite[Lemma 5.2]{surfaces}. Summing over all such $n_1,...,n_{l+1}$, $\overline{x}\in X_n$ and multiplying by $n!$ yields equation ~\eqref{eq:general-formula}. The interested reader may find a complete proof in \cite{thesis}.
\end{proof}

Let $\theta\in[0,1]$ and denote by $\o(\theta):=e^{2\pi i\theta}$. The motivation for the next theorem is explained in Problem 2.

\begin{thm}\label{thm:w-signature-basis}
Let $AP_n$ be the Abelian subgroup of $P_n$ generated by $\eta_{i,n}$, and denote by
$$V_{AP}:=\{\widetilde{\varphi}|_{AP_n}:AP_n\to\R\ \textrm{where}\ \widetilde{\varphi} \textrm{ is a  homogeneous quasi-morphism on $B_n$}\}.$$
Then
$$\mathcal{B}=\left\{\widetilde{sign}_{\o\left(\frac{1}{2}\right)},
2\widetilde{sign}_{\o\left(\frac{1}{2}\right)}
-3\widetilde{sign}_{\o\left(\frac{1}{3}\right)},
\ldots,(n-1)\widetilde{sign}_{\o\left(\frac{1}{n-1}\right)}-
n\widetilde{sign}_{\o\left(\frac{1}{n}\right)}\right\}$$
is a basis for $V_{AP}$.
\end{thm}

\begin{proof}
If $n=2$, then
$\widetilde{sign}_{\o\left(\frac{1}{2}\right)}(\eta_{2,2})=2$, and the proof follows. Let $n>2$.

\begin{lem}\label{lem:w-signature-of-eta-i-n}
Let $\theta\in[0,1]\cap\mathbb{Q}$. Then:
\begin{multline*}
\widetilde{sign}_{\o(\theta)}(\eta_{i,n})=\left\{
                                          \begin{array}{c}\begin{aligned}
                                            &4(i-1)\theta,  &\rm{if}&\ 0\leq\theta \leq\frac{1}{i},\\\\
                                            &4(l-1)(1-\theta),
                                            &\rm{if}&\
                                           \frac{l-1}{i}\leq\theta \leq\frac{l-1}{i-1}, \quad 2\leq l \leq i-1,\\\\
                                            &4(i-l)\theta,  &\rm{if}& \  \frac{l-1}{i-1}\leq\theta \leq\frac{l}{i}, \quad 2\leq l \leq i-1,\\\\
                                            &4(i-1)(1-\theta),  &\rm{if}& \ \frac{i-1}{i}\leq\theta \leq 1.\\
                                            \end{aligned}
                                          \end{array}
                                        \right.
\end{multline*}
\end{lem}
\begin{proof}
Recall that the torus link $K(p,q)$ has a braid representation $(\s_1\cdot\ldots\cdot\s_{p-1})^q$. Note that
$$\widetilde{sign}_{\o(\theta)}(K(n,n))=\widetilde{sign}_{\o(\theta)}(\eta_{2,n}\cdot\ldots\cdot\eta_{n,n})=
\sum\limits_{i=2}^n\widetilde{sign}_{\o(\theta)}(\eta_{i,n}).$$
Hence for each $i\leq n$ we have
\begin{equation}\label{eq:torus-link}
\widetilde{sign}_{\o(\theta)}(\eta_{i,n})=\widetilde{sign}_{\o(\theta)}(K(i,i))-\widetilde{sign}_{\o(\theta)}(K(i-1,i-1)).
\end{equation}

In \cite{w-signature} Gambaudo and Ghys proved the following
\begin{prop}[\cite{w-signature}, Proposition 5.2]\label{prop:w-sign-calc}
Let $\theta\in[0,1]\cap\mathbb{Q}$. Then
$$
\widetilde{sign}_{\o(\theta)}(\s_1\cdot\ldots\cdot\s_{n-1})=2(n-2l+1)\theta +\frac{2l(l-1)}{n}
$$
for $\frac{l-1}{n}\leq\theta\leq\frac{l}{n}$, where $1\leq l\leq n$.
\end{prop}
An immediate consequence of Proposition \ref{prop:w-sign-calc} is that for $1\leq l\leq n$ and $\frac{l-1}{n}\leq\theta\leq\frac{l}{n}$ we have
\begin{equation}\label{eq:sign-w}
\widetilde{sign}_{\o(\theta)}(K(n,n))=2n(n-2l+1)\theta +2l(l-1).
\end{equation}
A simple calculation combined with equalities \eqref{eq:torus-link} and \eqref{eq:sign-w} yields the proof of this lemma.
\end{proof}

Let us finish the proof of the theorem.
It follows from Lemma \ref{lem:w-signature-of-eta-i-n} that for each $3\leq i\leq n$ and for each $2\leq j\leq i-1$ we have
$$\left((i-1)\widetilde{sign}_{\o\left(\frac{1}{i-1}\right)}-
i\widetilde{sign}_{\o\left(\frac{1}{i}\right)}\right)(\eta_{i,i})=-4,$$
$$\left((i-1)\widetilde{sign}_{\o\left(\frac{1}{i-1}\right)}-
i\widetilde{sign}_{\o\left(\frac{1}{i}\right)}\right)(\eta_{j,i})=0.$$
Thus $\mathcal{B}$ is a basis for $V_{AP}$.
\end{proof}

\begin{pr*}
The following problem is motivated by a question posed to us by L.Polterovich. Denote by $\cD_{aut}\subset \cD$ the set of area-preserving diffeomorphisms generated by autonomous Hamiltonians. Observe that, by Theorems \ref{thm:general-formula} and \ref{thm:w-signature-basis}, one can find a linear combination of homogeneous quasi-morphisms
$\widetilde{sign}_\o$ on $B_n$ so that the corresponding homogeneous quasi-morphism $\widetilde{\Phi}: \cD\to\R$ vanishes on $\cD_{aut}$.
On the other hand, using a construction from \cite{surfaces} (completely different from the one described in this paper) one
can construct more homogeneous quasi-morphisms $\widetilde{\Psi}_\a$ on $\cD$ that vanish on $\cD_{aut}$. It would be interesting to check whether the quasi-morphism $\widetilde{\Phi}$ is non-trivial and, if so, whether it is a linear combination of the quasi-morphisms $\widetilde{\Psi}_\a$. So far we have not been able to compute or estimate the value of {\it any} of the Gambaudo-Ghys quasi-morphisms coming from quasi-morphisms on $B_n$ ($n>2$) on {\it any} area-preserving diffeomorphism that is not generated by a flow preserving a foliation.
\end{pr*}

\subsection{Asymptotic behavior of an induced signature quasi-morphism}

Recall that $\widetilde{sign}$ is a homogeneous
quasi-morphism on $B_n$. The induced family of homogeneous quasi-morphisms on $\cD$ is denoted by
$\widetilde{Sign}_{n,\D^2}$. In \cite{surfaces} Gambaudo and Ghys  proved that all
of them are non-trivial and linearly independent.

\begin{prop}\label{prop:signature-quasimorphism}
For each $n\geq 2$ we have
\begin{equation*}
\widetilde{Sign}_{n,\D^2}(h_1)=
n\pi^{n}\int\limits_T\left(1+4(n-1)J-(1-4J)^{n-1}\right)\hbar^\prime(J)dJ,
\end{equation*}
where $h_1$ is the time-one Hamiltonian flow defined by a Morse-type Hamiltonian $H$.
\end{prop}
\begin{proof}
In \cite{surfaces} Gambaudo and Ghys showed that
\begin{equation*}
\widetilde{sign}(\eta_{i,n})=\left\{
\begin{array}{c}\begin{aligned}
&i,\hspace{10.5mm} \rm{if}\hspace{2mm} \textit{i}\hspace{2mm} \rm{is}\hspace{2mm} \rm{even},\\
&i-1,\hspace{4.4mm} \rm{if}\hspace{2mm} \textit{i}\hspace{2mm} \rm{is}\hspace{2mm} \rm{odd}.\\
 \end{aligned}
 \end{array}
 \right.
\end{equation*}
A simple computation yields
\begin{align*}
&\sum\limits_{i=2}^n\widetilde{sign}_{n}(\eta_{i,n})i\dbinom{n}{i}J^{i-1}\left(\frac{1}{2}-J\right)^{n-i}=
\frac{n}{2}\left(\left(\frac{1}{2}\right)^{n-1}+4(n-1)J\left(\frac{1}{2}\right)^{n-1}\right)-\\
&\frac{n}{2}\left(\frac{1}{2}-2J\right)^{n-1}.
\end{align*}
Now the proof follows immediately from  Theorem \ref{thm:general-formula}.
\end{proof}

The following theorem is used in the proof of the main result in this subsection (Theorem \ref{thm:asymptotics}). The proof relies on the estimates from \cite[Theorem 1]{Gambaudo-Lagrange}. It is technical and may be found in \cite[Theorem 3.3.5]{thesis}.

\begin{thm}[\cite{thesis}]\label{thm:Morse-Hamiltonian}
Let $H\in\cH$ and $\{H_k\}_{k=1}^\infty$ be a sequence of functions such that each $H_k\in\cH$ and $H_k\xrightarrow[k\rightarrow\infty]{}H$ in $C^1$-topology. Let $h_1$ and $h_{1,k}$ be the time-one Hamiltonian flows generated by $H$ and $H_k$ respectively. Let $\widetilde{\varphi}_n:B_n\to\R$ be a homogeneous quasi-morphism, and $\widetilde{\Phi}$ the corresponding homogeneous quasi-morphism on $\cD$. Then
$$\lim\limits_{k\to\infty}\widetilde{\Phi}(h_{1,k})=\widetilde{\Phi}(h_1).$$
\end{thm}

Now we are ready to prove our main theorem in this subsection.

\begin{thm}\label{thm:asymptotics}
For each $h_1\in\cD$ generated by an autonomous
Hamiltonian $H$ we have
\begin{equation*}
\lim\limits_{n\to\infty}\frac{\widetilde{Sign}_{n,\D^{2}}(h_1)}{\pi^{n-1}n(n-1)}=\mathcal{C}(h_1),
\end{equation*}
where $\C$ is the Calabi homomorphism.
\end{thm}
\begin{proof}
\textbf{Step 1.} Suppose that  $h_1$ is generated by a  Morse-type
Hamiltonian $H$. It follows from Proposition \ref{prop:signature-quasimorphism} that
\begin{align*}
&\lim\limits_{n\to\infty}\frac{\widetilde{Sign}_{n,\D^2}(h_1)}{\pi^{n-1}n(n-1)}=
4\pi\int\limits_T J\hbar^\prime(J)dJ+\\
&\pi\lim\limits_{n\to\infty}\left(\int\limits_T \hbar^\prime(J)dJ/(n-1)-
\int\limits_T (1-4J)^{n-1}\hbar^\prime(J)dJ/(n-1)\right).
\end{align*}
The limit of the first integral in the equation above equals to zero, hence
\begin{equation*}
\lim\limits_{n\to\infty}\frac{\widetilde{Sign}_{n,\D^2}(h_1)}{\pi^{n-1}n(n-1)}=
4\pi\int\limits_T J\hbar^\prime(J)dJ-\lim\limits_{n\to\infty}\pi\left(\int\limits_T (1-4J)^{n-1}\hbar^\prime(J)dJ/(n-1)\right).
\end{equation*}
Note that $J$ is an action coordinate on each edge of $T$. This yields
$$0\leq J\leq\frac{\area(\D^2)}{2\pi}=\frac{\pi}{2\pi}=\frac{1}{2}.$$
It follows that $|1-4J|\leq 1$ and hence
$$0\leq\left|\lim\limits_{n\to\infty}\pi\left(\int\limits_T (1-4J)^{n-1}\hbar^\prime(J)dJ/(n-1)\right)\right|\leq
\lim\limits_{n\to\infty}\pi\left(\int\limits_T |\hbar^\prime(J)|dJ/(n-1)\right)=0.$$
Therefore
\begin{equation*}
\lim\limits_{n\to\infty}\frac{\widetilde{Sign}_{n,\D^2}(h_1)}{\pi^{n-1}n(n-1)}=
4\pi\int\limits_T J\hbar^\prime(J)dJ=-4\pi\int\limits_T \hbar(J)dJ=\C(h_1).
\end{equation*}
The second equality is just integration by parts, and the proof of the third equality is exactly the same as the proof of \cite[Proposition 2.2]{enlacement}.\\
\textbf{Step 2.} Suppose that  $h_1$ is generated by any $H\in\cH$. Then by \cite[Theorem 2.7]{Milnor} there exists a sequence $\{H_k\}_{k=1}^\infty$ of Morse-type functions, which converge in $C^1$-topology to $H$. Now the proof follows from Theorem \ref{thm:Morse-Hamiltonian}.
\end{proof}

\section*{Acknowledgments.} First of all, I would like to thank Michael
Entov, who has introduced me to this subject, guided and helped me a
lot while I was working on this paper. I would like to thank  Michael Polyak for his valuable suggestions
in the course of my work on this paper. I would like to thank  Leonid Polterovich for helpful comments and support.

Department of Mathematics, Vanderbilt University, Nashville, TN\\
\emph{E-mail address:} \verb"michael.brandenbursky@vanderbilt.edu"

\end{document}